\documentclass[conference, 11pt]{IEEEtran}
\IEEEoverridecommandlockouts

\usepackage{graphicx,color,placeins,float,tabularx,colortbl,amsmath,amssymb,epsfig,amsthm,cite}
\usepackage[pdfusetitle, pdfauthor={Tarek Aziz Lahlou, MIT}]{hyperref}
\hypersetup{colorlinks = true,linkcolor=black,citecolor=blue,urlcolor=blue,anchorcolor=blue}

\begin{document}  \onecolumn 

\title{Trading Accuracy for Numerical Stability: Orthogonalization, Biorthogonalization and Regularization}

\author{\IEEEauthorblockN{Tarek A. Lahlou and Alan V. Oppenheim\thanks{The authors wish to thank Analog Devices, Bose Corporation, and Texas Instruments for their support of innovative research at MIT and within the Digital Signal Processing Group.}} \IEEEauthorblockA{Digital Signal Processing Group  \\ Massachusetts Institute of Technology }}

\maketitle

\begin{abstract}
	This paper presents two novel regularization methods motivated in part by the geometric significance of biorthogonal bases in signal processing applications. These methods, in particular, draw upon the structural relevance of orthogonality and biorthogonality principles and are presented from the perspectives of signal processing, convex programming, continuation methods and nonlinear projection operators. Each method is specifically endowed with either a homotopy or tuning parameter to facilitate tradeoff analysis between accuracy and numerical stability. An example involving a basis comprised of real exponential signals illustrates the utility of the proposed methods on an ill-conditioned inverse problem and the results are compared to standard regularization techniques from the signal processing literature. 
\end{abstract}

\tableofcontents

\IEEEpeerreviewmaketitle

\newpage
\section{Introduction}
 	Signal processing algorithms, notably those in machine learning and big data applications, make ubiquitous use of dense and sparse linear algebra subroutines \cite{SPBigData, DeepLearning}. Numerical errors often arise while implementing such subroutines and are typically explained using perturbation analysis tools. For example, finite-precision effects such as coefficient quantization are quantifiable by invoking sensitivity theorems after organizing the computation into a linear or nonlinear signal-flow graph, as discussed in, e.g.,~\cite{os}. Also, the magnification and compounding of propagated errors are commonly understood via the condition number of a functional representation of the processing system. These unavoidable sources of error contribute significantly to the total accumulated error in an obtained solution, especially as problem sizes continue to scale. Algorithms specifically designed to use orthogonality principles alleviate these issues and, among other reasons, receive extensive use in practice. Moreover, orthogonality in other contexts such as quantum measurement and asynchronous distributed optimization continue to provide inspiration for new signal processing algorithms, architectures, and applications \cite{yoninaQuantum,BLOpt}.

 	Signal processing techniques exploiting natural and efficient representations for signal models with inherent geometric structure often use non-orthogonal bases and frames that may generally suffer from conditioning issues, e.g.~finite rate of innovation \cite{FRI}, compressed sensing \cite{cs}, and transient signal analysis \cite{Prony6}-\cite{prony2}. 
 	% \cite{Prony6}\cite{LahlouOppenheimGlobalSIP2014}\cite{ExponentialBook}\cite{prony2}
 	To address these issues, we restrict our attention in this paper to two classical approaches: orthogonalization and regularization. Traditional orthogonalization routines typically destroy the topological relevance of an original basis to the problem at hand while simultaneously weakening the informational significance of the solution obtained. The orthogonal Procrustes problem has been used for subspace whitening and machine learning in the matrix setting, i.e.~generating the unitary matrix that optimally maps one matrix to another in the Frobenius norm sense\cite{whiten1,whiten2}. In this paper, we specifically translate and explicate these results using bases and further extend them as they pertain to biorthogonalization, nonlinear projections, and novel regularization methods. The proposed regularization methods permit natural characterizations of accuracy and stability enabling standard tradeoff analysis to identify an acceptable balance. Regularization methods in signal processing have historically focused on directly obtaining a solution biased toward prespecified signal properties. In this paper, however, we emphasize the effect the proposed methods have on the linear problems Euclidean structure itself. 

	This paper is outlined as follows: In Section~\ref{sec:topology} we collect together geometric facts surrounding orthogonality and biorthogonality that will be referenced throughout. The orthogonalization and biorthogonalization of a linear system via a nonlinear projection operator is the central focus of Section~\ref{sec:orth-biorth}. Two regularization methods for balancing numerical robustness in the sense of relative conditioning and accuracy in the sense of geometric structure are developed in Section~\ref{sec:acc-sta}. The proposed methods are then applied to the regularization of an ill-conditioned basis composed of real exponential signals in Section~\ref{sec:exmpl} as a numerical example. Throughout the presentation, we draw no distinction between various linear problems and the mathematical relationships connecting them, e.g.~matched filtering, solving linear systems, changing bases, etc. In this way, the insights and methods presented may be readily applied to any linear problem equally well.

	\subsection{Notational conventions}\label{subsec:notation}
		Vectors are denoted using underscores with subscripts indexing vectors as opposed to entries, i.e.~$\underline{a}_1$ and $\underline{a}_2$ represent two distinct vectors. Boldface letters denote a \emph{system} or collection of $N$ vectors, typically constituting a basis for $\mathbb{R}^N$, i.e.~$\mathbf{a} = \{\underline{a}_{k}\}_{k=1}^{N}$. A biorthogonal set is denoted using $\sim$, i.e.~$\widetilde{\mathbf{a}}$ is biorthogonal to $\mathbf{a}$. The $p$-norm of a vector $\underline{a}$ for $p\geq1$ is denoted $\| \underline{a}\|_{p}$. A capital letter denotes the matrix form of a set of vectors, i.e.~$A$ has columns $\mathbf{a}$. The Kronecker delta is denoted $\delta_{k,j}$. For clarity we assume all systems are real-valued bases; rigorous extensions to complex-valued bases and frames follow analogously, e.g. by using appropriate conjugation and restricting arguments to the range of a system. 

\section{The topology of the Stiefel manifold and $O(N)$} \label{sec:topology}
	Let $\mathbf{a}$ and $\mathbf{b}$ denote two general systems and define $\epsilon_{BO}$ as the nonnegative measure of their biorthogonality given by
	\begin{eqnarray} 
		\epsilon_{BO}\left(\mathbf{a}, \mathbf{b} \right) & \triangleq & \sum_{k}\sum_{j} \left( \left\langle \underline{a}_{k}, \underline{b}_{j} \right\rangle - \delta_{k,j} \right)^{2} \label{eq:biorthogonality-measure}
	\end{eqnarray}
	which achieves a unique global minimum of zero provided that $\mathbf{b}=\widetilde{\mathbf{a}}$. It is then straightforward to conclude that
	\begin{equation}
		\epsilon_{BO}\left(\mathbf{a}, \mathbf{a} \right)  = 0 \iff \mathbf{a} \mbox{ is an orthonormal system}.
	\end{equation}
	The Stiefel manifold $V_{k,N}$ is the orthogonal groups principal homogenous space $(k\!=\!N)$ and is written using \eqref{eq:biorthogonality-measure} as \cite{topology}
	\begin{equation} \label{eq:ortho-bases}
		\displaystyle{V_{N,N} = \left\{\mathbf{h} \in \mathbb{R}^{N} \colon \epsilon_{BO}(\mathbf{h},\mathbf{h})=0\right\}.} 
	\end{equation}
	Represented as a quadratic form, \eqref{eq:biorthogonality-measure} is convex quadratic with positive signature, i.e. all eigenvalues are positive valued, thus any optimization problem that uses \eqref{eq:biorthogonality-measure} as a cost function or \eqref{eq:ortho-bases} as a feasible set is not precluded from being a convex problem.
		
	The orthogonal group of dimension $N$, denoted by $O(N)$ and described in set-builder notation using the matrix form of $\mathbf{h}$ as
	\begin{equation}
		O(N) = \left\{ H\in\mathbb{R}^{N \times N} \colon H^TH = HH^T = I_N \right\}
	\end{equation}
	where $I_N$ is the identity operator on $\mathbb{R}^{N}$, consists of two connected components $\mathcal{M}_{1}$ and $\mathcal{M}_{-1}$. Furthermore, $\mathcal{M}_{1}$ and $\mathcal{M}_{-1}$ each form a disjoint smooth manifold in $\mathbb{R}^{N \times N}$ characterized by
	\begin{equation}
		\mathcal{M}_{i} = \left\{H\in O(N) \colon \text{det}(H) = i\right\}, \hspace{1em} i = \pm 1 \label{eq:manifolds}
	\end{equation}
	where $\text{det}(\cdot)$ denotes the determinant operator. Note that $\mathcal{M}_{1}$ is precisely the special orthogonal group $SO(N)$. From \eqref{eq:manifolds} and the definition of a connected component it immediately follows that
	\begin{equation}
		H\in \mathcal{M}_{i}    \Leftrightarrow   \widetilde{H} \in \mathcal{M}_{i}, \hspace{1em} i = \pm 1. \label{eq:whichManifold}
	\end{equation}
	We conclude this section by noting that appropriately selecting a base point for $V_{N,N}$ establishes a one-to-one correspondence with $O(N)$ and, therefore, justifies interchanging between the two notational descriptions of an orthogonal system, i.e.~$\mathbf{h}$ and $H$.

\section{Orthogonalization and biorthogonalization}\label{sec:orth-biorth}
	Let $\mathbf{a}$ and $\mathbf{b}$ continue to denote two general systems and define $\epsilon_{LS}$ as the nonnegative Euclidean measure of their  distance given by
	\begin{eqnarray}
		\epsilon_{LS}\left(\mathbf{a}, \mathbf{b} \right) & \triangleq & \sum_{k} \|\underline{a}_{k} - \underline{b}_{k} \|_{2}^{2} \label{eq:ls-measure}\\ 
		& = & \sum_{k} \left\langle \underline{a}_{k} - \underline{b}_{k}, \underline{a}_{k} - \underline{b}_{k} \right\rangle 
	\end{eqnarray}
	which achieves a unique global minimum of zero provided that $\mathbf{b} = \mathbf{a}$. In the sequel we shall denote by $\mathbf{g}$ a general prespecified system.

	\subsection{Optimal orthogonalization problems}
		Consider the optimization problem that identifies the nearest orthogonal system to $\mathbf{g}$ in the sense of \eqref{eq:ls-measure}, or equivalently the least-squares orthogonalization of $\mathbf{g}$, given by
		\begin{equation}
			\widehat{\mathbf{h}}_{\mathbf{1}} = \arg \min_{\mathbf{h}} \epsilon_{LS}\left(\mathbf{g}, \mathbf{h}\right) \mbox{ s.t. } \epsilon_{BO}\left(\mathbf{h}, \mathbf{h} \right) = 0.            				\tag{P1}\label{eq:P1}
		\end{equation} 
		Similarly, consider the optimization problem that identifies the nearest orthogonal system to the biorthogonal system $\widetilde{\mathbf{g}}$ in the sense of \eqref{eq:ls-measure} given by 
		\begin{equation}
			\widehat{\mathbf{h}}_{\mathbf{2}} = \arg \min_{\mathbf{h}} \epsilon_{LS}\left(\widetilde{\mathbf{g}}, \mathbf{h}\right) \mbox{ s.t. } \epsilon_{BO}\left(\mathbf{h}, \mathbf{h} \right) = 0.             \tag{P2}\label{eq:P2}
		\end{equation} 
		We collectively refer to \eqref{eq:P1} and \eqref{eq:P2} as optimal orthogonalization problems and further interpret their solutions as orthogonal systems that maximally preserve the geometric structure inherent to $\mathbf{g}$ and $\widetilde{\mathbf{g}}$, respectively.

	\subsection{Optimal biorthogonalization problems}
		Consider the optimization problem that identifies the orthogonal system which is maximally dual to $\mathbf{g}$ in the sense of \eqref{eq:biorthogonality-measure} given by
		\begin{equation}
			\widehat{\mathbf{h}}_{\mathbf{3}} = \arg \min_{\mathbf{h}} \epsilon_{BO}\left(\mathbf{g}, \mathbf{h}\right) \mbox{ s.t. } \epsilon_{BO}\left(\mathbf{h}, \mathbf{h} \right) = 0.            				\tag{P3}\label{eq:P3}
		\end{equation}
		Finally, consider the optimization problem that identifies the orthogonal system which is maximally dual to the biorthogonal system $\widetilde{\mathbf{g}}$ in the sense of \eqref{eq:biorthogonality-measure} given by
		\begin{equation}
			\widehat{\mathbf{h}}_{\mathbf{4}} = \arg \min_{\mathbf{h}} \epsilon_{BO}\left(\widetilde{\mathbf{g}}, \mathbf{h}\right) \mbox{ s.t. } \epsilon_{BO}\left(\mathbf{h}, \mathbf{h} \right) = 0.             \tag{P4}\label{eq:P4}
		\end{equation} 
		We collectively refer to \eqref{eq:P3} and \eqref{eq:P4} as optimal biorthogonalization problems and further interpret their solutions as the orthogonal systems that, switching to the matrix form notation introduced in Subsection~\ref{subsec:notation}, maximally invert $G$ and $\widetilde{G}$, respectively.

		\begin{figure}[t]
		 	\centering
		  	\centerline{\includegraphics[width=9cm]{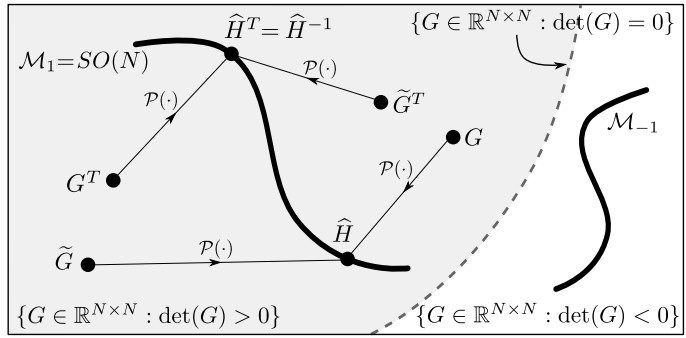}}
		  	\caption{An illustration of the relationships between the matrix forms of $\mathbf{g}$ and $\widetilde{\mathbf{g}}$ and their adjoints and their projections onto the orthogonal group via \eqref{eq:soln} where the projection operator is denoted using $\mathcal{P}(\cdot)$. The system $\mathbf{g}$ is chosen such that $\text{det}(G) > 0$ and therefore the projections are onto $\mathcal{M}_{1}$.}
			\label{fig:manifold-projections}
		\end{figure}

	\subsection{Equivalence of $\epsilon_{LS}$ and $\epsilon_{BO}$ over $V_{N,N}$}\label{subsec:equivalence}
		The measures \eqref{eq:biorthogonality-measure} and \eqref{eq:ls-measure}, respectively interpreted as optimal biorthogonalization and orthgonalization cost functions, produce generally unequal cost surfaces over the space of all $N$-dimensional bases. However, for fixed $\mathbf{g}$ and orthogonal $\mathbf{h}$, it follows that
		\begin{eqnarray}
			\epsilon_{BO}(\mathbf{g}, \mathbf{h})	& = & 	\sum_{k,\,j}\left(\langle \underline{g}_{k},\underline{h}_{j}\rangle \left[\langle \underline{g}_{k},\underline{h}_{j}\rangle -2\delta_{k,j}\right]+\delta_{k,j}^{2}\right) \label{eq:proof1} \\ 
			& = & \sum_{k}	\left(\langle \underline{g}_{k},\underline{g}_{k}\rangle -2\langle \underline{g}_{k},\underline{h}_{k}\rangle +\langle \underline{h}_{k}, \underline{h}_{k} \rangle\right) \label{eq:proof2}\\
			& = & \sum_{k}\langle \underline{g}_{k}-\underline{h}_{k},\underline{g}_{k}-\underline{h}_{k}\rangle   \label{eq:proof3}\\
			& = & \epsilon_{LS}(\mathbf{g}, \mathbf{h})
		\end{eqnarray}
		where the restriction of $\mathbf{h}$ to $V_{N,N}$ is used to obtain \eqref{eq:proof2} from \eqref{eq:proof1}, in particular since $\langle\underline{a}, \underline{b}\rangle = \sum_{j} \langle\underline{a}, \underline{h}_{j}\rangle\langle\underline{h}_{j}, \underline{b}\rangle$ and $\delta_{k,j}=\langle\underline{h}_{k}, \underline{h}_{j}\rangle$.

	\subsection{Analytic solution and projection interpretation of \eqref{eq:P1}-\eqref{eq:P4}}
		In Subsection~\ref{subsec:equivalence}, we established that the solution to \eqref{eq:P1} solves \eqref{eq:P3} and likewise that the solution to \eqref{eq:P2} solves \eqref{eq:P4}. In this subsection, we first state and interpret the analytic expression that minimizes \eqref{eq:P1} followed by discussing one way \eqref{eq:P1} and \eqref{eq:P2} may be argued to have the same minimizer. Indeed, the matrix form of the orthogonal system $\widehat{\mathbf{h}}_{\mathbf{1}}$ which solves \eqref{eq:P1} is given~by 
		\begin{equation}
			\widehat{H}_{1} = UV^{T}\label{eq:soln}
		\end{equation}
		where $G=U\Sigma V^{T}$ is the Singular Value Decomposition (SVD) of the matrix form of $\mathbf{g}$. Observe that verifying the orthogonality of \eqref{eq:soln} folllows immediately from the closure property of $O(N)$; a sketch of the full argument that \eqref{eq:soln} minimizes \eqref{eq:P1} is deferred to the Appendix. Interpreting \eqref{eq:soln} as a nonlinear operator acting upon $G$, it is straightforward to verify idempotency, hence \eqref{eq:soln} is the projection of $G$ onto $O(N)$. Figure~\ref{fig:manifold-projections} illustrates this, i.e.~that the action of \eqref{eq:soln} on $G$, $\widetilde{G}$ and their transposes is a projection onto the orthogonal group. From the perspective of bases, generating \eqref{eq:soln} consists of $N$ radial projections along the semiaxes of the ellipsoid described by the SVD of $G$ such that the associated singular values become unity. 

		An equivalence between \eqref{eq:P1} and \eqref{eq:P2} is established by using the matrix and basis representation equivalence discussed earlier where the matrix form of $\widetilde{\mathbf{g}}$, i.e.~the linear functionals biorthogonal to $\mathbf{g}$ in the dual linear space, corresponds to the inverse of the adjoint of $G$. The desired result is obtained by straightforward manipulations that further simplify \eqref{eq:P2}, in particular by exploiting SVD properties. For arbitrary systems $\mathbf{g}$, we comment that it is possible to identify \emph{a priori} which manifold \eqref{eq:soln} belongs to, i.e.~combining the manifold characterizations in \eqref{eq:manifolds} and using the connected component property in \eqref{eq:whichManifold} we conclude that
		\begin{equation}
			\begin{array}{lccc}
				\widehat{H}\in \mathcal{M}_{i}  	& 	\iff 	& 	\text{sgn}\left(\text{det}(G)\right) = i, & i = \pm1.    \\ 
			\end{array}
		\end{equation}

	\subsection{Incremental projection methods for nearly orthogonal $\mathbf{g}$}
		Although many numerical algorithms exist to compute \eqref{eq:soln}, we call special attention to the case of $G$ being ``nearly'' orthogonal. Let $R = G^TG - I_N$ be a Grammian residual. Re-writing \eqref{eq:soln} as 
		\begin{equation}
			\widehat{H} = G(I_N + R)^{-\frac{1}{2}},
		\end{equation}
		we proceed assuming the spectral norm of $R$ is less than $1$ which may be efficiently certified using, e.g.~Gershgorin's circle theorem~\cite{circleThm}. After some straightforward manipulations we obtain the power series
		\begin{equation}
			\widehat{H} = G\left[\sum_{n=0}^{\infty}{-\frac{1}{2} \choose n}R^{n}\right] \label{eq:seriesExansion}
		\end{equation}
	 	and thus the solution to \eqref{eq:P1}-\eqref{eq:P4} may be generated to any desired level of accuracy without explicitly computing an SVD by truncating \eqref{eq:seriesExansion} to an appropriate number of terms.

\section{Trading accuracy for numerical stability} \label{sec:acc-sta} 
	Many applications of signal processing as well as engineering more broadly require considerably less accuracy than what is available on modern computational platforms. By pairing systems in $V_{N,N}$ with backward stable algorithms, computation can be made to not magnify propagated errors. Motivated by these observations, we next propose two regularization methods that facilitate accuracy and stability tradeoffs readily explained through homotopy maps and convex analysis, respectively. We comment upfront that while our treatment focuses on interpretation, the utility of these methods is in regularizing ill-conditioned problems while preserving some of the original systems Euclidean structure as is illustrated in Section~\ref{sec:exmpl}.

	Regularization methods are often described initially using the language of optimization theory; the methods surveyed in Subsection~\ref{sec:methods} serve as examples. Some regularization techniques, primarily those related to least squares formulations, admit solutions in closed-form. Setting these special cases aside, the remaining techniques conventionally solve directly for a solution to the problem at hand without first generating a regularized system. The methods proposed in this section specifically make this intermediary step, i.e.~they explicitly generate a regularized linear system from which a solution is then generated, potentially making further use of any regularization techniques belonging to the class just mentioned in which the regularized system takes the place of the original system.

	\subsection{Homotopic continuation formulation}\label{sec:oper}
		Motivated by the extensive use of floating-point arithmetic in big data computation, we measure the stability of a system $\mathbf{g}$ through its relative condition number, i.e.~the ratio of extremal singular values of $G$. Loosely speaking, when this measure is on the order $10^{k}$ then only  $d-k$  digits of the computed solution are reliable where $d$ is the maximum number available on the chosen computational platform \cite{CONDNUMDIGITS}. This loss in accuracy is inherent to the system itself and says nothing about further degradations that accumulate due to a particular algorithms stability nor to errors that arise from round-off noise and coefficient quantization. 

		We first introduce a continuation scheme to generate a parameterized family of regularized systems along the trajectory corresponding to the smooth deformation of $\mathbf{g}$ into $\widehat{\mathbf{h}}_{\mathbf 1}$. The utility of this scheme for solving a linear system of equations $G\underline{x} = \underline{y}$, with known $\underline{y}$, is to sequentially or approximately solve for  $\underline{x}$ while substituting different systems along this continuum for $G$. The parameterization of the system mitigates numerical conditioning issues at the expense of accuracy, i.e.~the continuation terminates when a system associated with an acceptable solution is identified. The function describing these deformations is referred to as a \emph{homotopy map}. In this paper, we use the specific homotopy map
		\begin{equation}
			\widehat{\mathbf{h}}_{\mathbf{5}}(\rho, \mathbf{g}) \!\triangleq\! \left\{ \underline{h}_{k} \colon \underline{h}_{k} = (1\!-\!\rho)\underline{g}_{k} + \rho\underline{z}_{k} \!\text{ where } \mathbf{z} \text{ solves (P1)}\right\} \label{eq:homotopy}
		\end{equation}
		where the parameter $\rho\in[0,1]$ balances the previously addressed tradeoff. The family of linear systems achieved using the homotopy map $\widehat{\mathbf{h}}_{\mathbf{5}}$ is  $\widehat{\mathbf{h}}_{\mathbf{5}}(\rho, \mathbf{g}) \leftrightarrow \widehat{H}_{5}$ for $\rho\in[0,1]$ and the standard implementation strategy is to track the solution to the problem using $\widehat{H}_{5}$ instead of $G$ starting from $(\rho,\mathbf{g})=(0,\mathbf{g})$ as $\rho$ progresses from $0$ to $\rho^\star$ where $\rho^\star$ generates the system associated with an acceptable solution. In Section~\ref{sec:exmpl}, solutions to a similar problem are obtained by discretizing $[0,1]$ and evaluating a performance metric for all values of $\rho$ in this interval. Standard line search methods help to efficiently identify the optimal homotopy value $\rho^\star$. 
		
		Since the condition number of $G$ is equal to the condition number of $\widetilde{G}$, and with the equivalence between the minimizer of \eqref{eq:P1} and \eqref{eq:P2}, we conclude the discussion in this subsection by noting that replacing  $\mathbf{g}$ with $\widetilde{\mathbf{g}}$ results in the same presentation with the geometric structure of the biorthogonal system $\widetilde{\mathbf{g}}$ being preserved instead.

	\subsection{Unconstrained optimization formulation}
		Prompted by the exposition on characterizing numerical accuracy and stability in the previous subsection, we next formulate an optimization problem whose cost function is the sum of the metrics \eqref{eq:ls-measure} and \eqref{eq:biorthogonality-measure}, i.e.~each summand is respectively minimized by $\widehat{\mathbf{h}}_{\mathbf{5}}(0, \mathbf{g})$ and $\widehat{\mathbf{h}}_{\mathbf{5}}(1,\mathbf{g})$. We formally write this problem as an unconstrained quartic convex optimization problem of the form
		\begin{equation}
			\widehat{\mathbf{h}}_{\mathbf{6}}(\rho,\mathbf{g}) \triangleq \arg \min_{\mathbf{h}} \epsilon_{LS}( \mathbf{g}, \mathbf{h}) +\rho\epsilon_{BO}(\mathbf{h}, \mathbf{h})             \label{eq:P5}
		\end{equation} 
		where $\rho$ again determines the previously addressed tradeoff. Indeed, consistent with the discussion surrounding the positive definite nature of \eqref{eq:biorthogonality-measure} represented as a quadratic form and the fact that the sum of unconstrained convex functions results in a convex function, it is straightforward to verify that any solution to  \eqref{eq:P5} that is locally optimal is also globally so. Although the cost function of \eqref{eq:P5} is purposefully constructed using the endpoints of the continuum of regularized systems achievable using \eqref{eq:homotopy}, it is readily possible to generate solutions to \eqref{eq:P5} that do not lie on this continuum nor are obtainable using \eqref{eq:homotopy} for any value of the homotopy parameter $\rho$.

		In response to the global optimality guarantees established above and for the sake of completeness, we conclude this section by including gradient primitives that may be used by any number of gradient-based nonlinear programming algorithms in solving \eqref{eq:P5}:
		\begin{eqnarray} 
			\nabla \epsilon_{LS}( \mathbf{g}, \mathbf{h}) \mid_{{\underline{h}}_{k}}\!\!\!\!\!\! & = & \!\!\!\!\! 2\sum_{j}(\langle \underline{g}_{j},\,\underline{h}_{k}\rangle -\delta_{j,k})\underline{g}_{j} \\
			\nabla \epsilon_{BO}(\mathbf{h}, \mathbf{h} ) \mid_{{\underline{h}}_{k}}\!\!\!\!\!\! & = & \!\!\!\!\! 4\left(\left\langle \underline{h}_{k},\underline{h}_{k}\right\rangle -1\right)\underline{h}_{k} + 2\sum_{j\neq k}\left\langle \underline{h}_{j},\underline{h}_{k}\right\rangle\underline{h}_{j}.
		\end{eqnarray}
		Note that an analytic expression for the system $\widehat{\mathbf{h}}_{\mathbf{6}}$ which minimizes \eqref{eq:P5} is currently unknown due to the cubic nature of the first order optimality conditions, i.e.~the nonlinear system of equations defined by setting the appropriately weighted sum of the two terms above equal to zero for each value of $k$.

\section{Numerical example}\label{sec:exmpl}
	Solving linear equations using regularization, a common approach to ill-posed problems, modifies the problem so that (i) the new problem is biased toward expected solutions and (ii) the previously mentioned numerical issues are reduced. In signal processing, regularization appears in many forms: Tikhonov regularization and Wiener filtering, total variation denoising in image processing, and basis pursuit denoising in compressive sensing. In this section, we use the methods developed in Section~\ref{sec:acc-sta} to solve a numerical example. 

	\subsection{Problem formulation}
		Let $E$ denote the matrix form of a real exponential basis $\mathbf{e}$, i.e.~$E$ is Vandermonde and, for $0 < \sigma_{1} < \cdots < \sigma_{N} < 1$, generated by
		\begin{equation} 
			\hspace{2em} E_{i,\,j} = \sigma_j^{i-1}, \hspace{1em}1 \leq i,j \leq N. \label{eq:VandermondeBasis}
		\end{equation}
		The relative condition number of $E$, and consequently $\widetilde{E}$ too, grows exponentially in $N$ \cite{VanderEst}. This fact justifies using regularization to solve a linear problem of the form 
		\begin{equation} \label{eq:exmpl} 
			E\underline{x} = \underline{y}
		\end{equation}
		for $\underline{x}$ where $\underline{y}$ is analytically generated according to a synthetic solution $\bar{\underline{x}}$ to help limit observed numerical effects to the regularization effect each method has on the matrix $E$ in \eqref{eq:exmpl}.  

	\subsection{Numerical methods for comparison}\label{sec:methods}
		This subsection briefly reviews three regularization methods from the signal processing literature, two of which originate in the compressive sensing setting, to compare against \eqref{eq:homotopy} and \eqref{eq:P5} in solving \eqref{eq:exmpl} \cite{cs}. The first method, Tikhonov regularization, is formulated as
		\begin{equation}\label{eq:Tikh}
			\underline{x}^{\star} = \arg \min_{\underline{x}} \|E\underline{x} - \underline{y}\|_2^2 + \|T\underline{x}\|_2^2. %pick $T = \rho I_N$
		\end{equation}
		where we proceed to take the Tikhonov matrix $T=\rho I_N$. The objective in \eqref{eq:Tikh} is convex quadratic, thus the solution has an analytic expression, namely $\underline{x}^{\star} = (E^TE+\rho^2I)^{-1}E^T\underline{y}$. Next, the basis pursuit denoising problem is formulated as 
		\begin{equation}\label{eq:bpdn}
			\underline{x}^{\star} = \arg \min_{\underline{x}} \|E\underline{x} - \underline{y}\|_2 + \rho \|\underline{x}\|_1 
		\end{equation}
		where $\rho$ balances the absolute size of the solution with the desired agreement of \eqref{eq:exmpl}. Finally, the Dantzig selector is given by 
		\begin{equation}\label{eq:ds}
			\underline{x}^{\star} = \arg \min_{\underline{x}} \|E^T(E\underline{x} - \underline{y})\|_\infty + \rho \|\underline{x}\|_1
		\end{equation}
		and may be solved using standard linear programming techniques. 

	\subsection{Numerical results}
		Our experimental setup is as follows: We generate $E$ on each trial by selecting $N=18$ values $\sigma_k$ uniformly at random from the interval $[0.1,0.9]$ and solve \eqref{eq:exmpl} using the various methods discussed. We select $\rho$ to within $10^{-6}$ for each method per trial to minimize the residual $\|\underline{x}^{\star}-\bar{\underline{x}}\|_2$. The solution $x^{\star}$ to \eqref{eq:homotopy} and \eqref{eq:P5} is found by directly solving the regularized equations via Gaussian elimination. Table~1 lists the average residual and value of $\rho$ taken over $10^7$ trials. Solving  \eqref{eq:exmpl} directly yields an average error of $8.64$. The condition number of $E$ for this experiment is lower bounded by $10^{16}$. In summary, we remark that \eqref{eq:homotopy} outperforms the alternatives and \eqref{eq:P5} outperforms all but \eqref{eq:Tikh} for the Vandermonde linear systems tested. The larger error for \eqref{eq:ds} is due in part to the condition number of $E^TE$ being the square of the condition number of $E$.

\begin{figure}[t]
	 	\centering
	  	\centerline{\includegraphics[width=12cm]{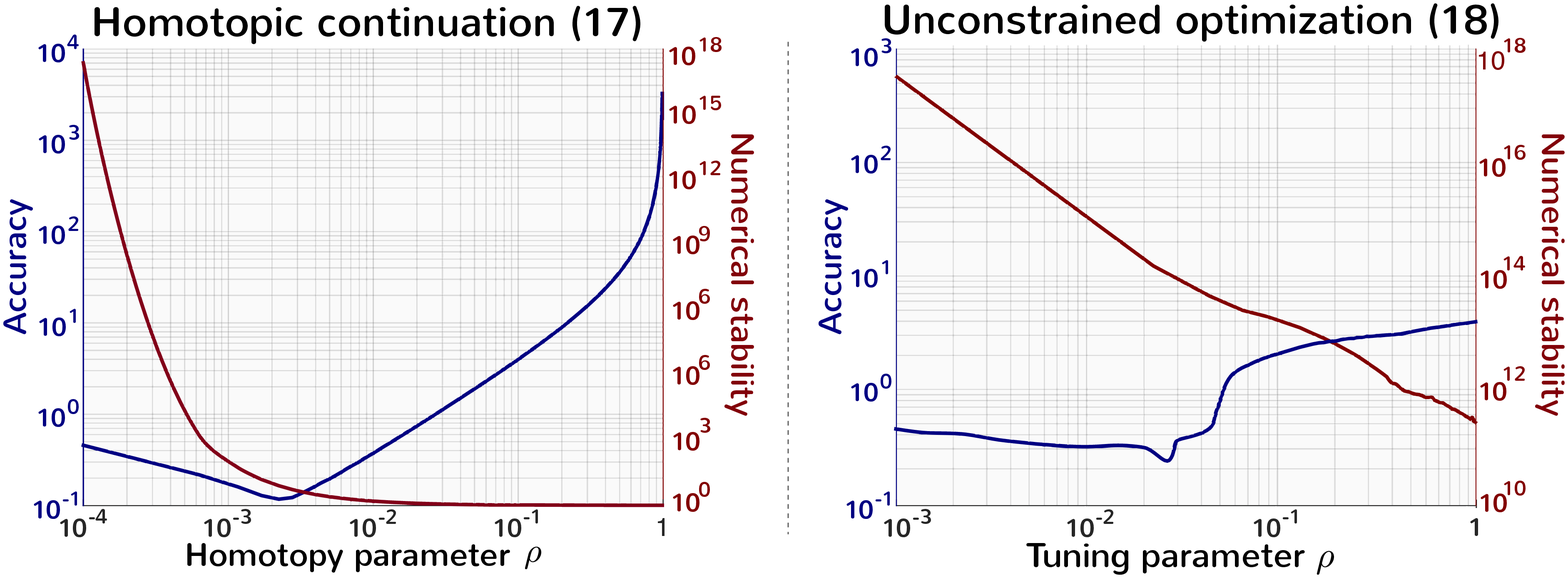}} %caption for fig 2
	  	\caption{Accuracy (residual error) versus numerical stability (relative condition number) tradeoff curves for \eqref{eq:homotopy} and \eqref{eq:P5}.}% figure below
		\label{fig:num-exmpl2} 
	\end{figure}

\begin{figure}[t]
	 	\centering
	  	\centerline{\includegraphics[width=12cm]{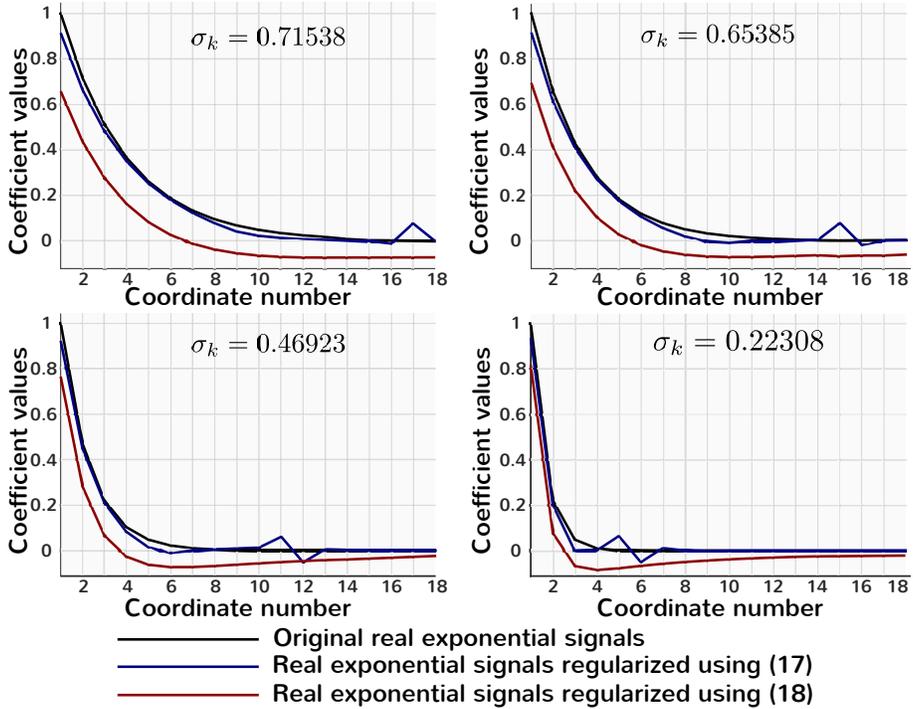}} % caption for fig 3
	  	\caption{An illustration of a subset of the real exponential basis elements $\mathbf{e}$ and the corresponding elements from the proposed regularization methods. } % text below
		\label{fig:num-exmpl} 
	\end{figure}

		Figure~\ref{fig:num-exmpl2} illustrates an accuracy versus numerical stability tradeoff curve for the proposed regularization methods where accuracy is given by the residual $\|\underline{x}^{\star}-\bar{\underline{x}}\|_2$ and numerical stability is given by the relative condition number. For both methods, it is clear that a non-zero value of $\rho$ is optimal corresponding to a system whose conditioning is several orders of magnitude less than $E$'s. 
		Figure~\ref{fig:num-exmpl} depicts a subset of regularized signals for the optimal selection of $\rho$ in Figure~\ref{fig:num-exmpl2}. As depicted and consistent with elements not shown, the regularized signals have similar structure to $\mathbf{e}$, i.e.~they decay at similar rates to their exponential counterparts, also depicted for comparison. Signals produced by \eqref{eq:P5}  retain a smooth nature with an additive offset and slight upward curvature and those produced by \eqref{eq:homotopy} contain small fluctuations that seem to alleviate conditioning~issues.


	\begin{center}
		\textsc{\small Table 1. Numerical comparison of regularization methods}
		\begin{tabular}{ c | c | c | c | c | c }
			Solution reference & \eqref{eq:homotopy} & \eqref{eq:P5} & \eqref{eq:Tikh} & \eqref{eq:bpdn} & \eqref{eq:ds} \\
			\hline
			\hline
			average $\|\underline{x}^{\star} -\bar{\underline{x}}\|_2$ 	& $0.096$ & $0.153$  &	$0.139$ & $0.147$ & $13.46$ \\
			average optimal $\rho$ 										& $0.002$ & $0.028$  &	$0.025$ & $0.015$ & $0.044$
		\end{tabular}
	\end{center}

\newpage
\section*{Appendix}
	Sketch of the argument that \eqref{eq:soln} solves \eqref{eq:P1}:\vspace{-0.05in}
	\begin{eqnarray*}
		\arg \min_{\mathbf{h} \in V_{N,N}} \epsilon_{LS}\left(\mathbf{g}, \mathbf{h}\right)  & = & \arg \min_{\mathbf{h} \in V_{N,N}} \sum_{k} \langle \underline{g}_k - \underline{h}_k, \underline{g}_k - \underline{h}_k \rangle \\
		 & = & \min_{\mathbf{h} \in V_{N,N}} \sum_{k} \left( \langle \underline{g}_k , \underline{g}_k \rangle + \langle \underline{h}_k, \underline{h}_k \rangle - 2 \langle \underline{g}_k, \underline{h}_k \rangle \right) \\
		 & = & \arg \max_{\mathbf{h}\in V_{N,N}}  \langle \underline{g}_{1}, \underline{h}_{1} \rangle  + \cdots + \langle \underline{g}_{N}, \underline{h}_{N} \rangle  
	\end{eqnarray*}
	\vspace{-0.15in} \\
	Switching to matrix notation for convenience, utilizing properties of a trace, and writing $G$ using its SVD $G=U\Sigma V^{T}$,  we proceed as:\vspace{-0.05in}
	\begin{eqnarray*}
		\arg \max_{\mathbf{h}\in V_{N,N}} \sum_{k} \langle \underline{g}_{k}, \underline{h}_{k} \rangle  & \rightarrow & \arg \max_{H\in O(N)} \text{trace}\left(V\Sigma U^T H\right) \\
		& = &  U\left( \arg \max_{H\in O(N)} \text{trace}\left(\Sigma H\right)\right)V^T \\
		& = & UV^T
	\end{eqnarray*}
\newpage
\FloatBarrier
\bibliographystyle{IEEEbib}
\bibliography{refs} %\nocite{*}

\end{document}